\documentstyle [12pt,amsfonts] {article}
\input epsf

\newcommand{\beq}{\begin{equation}}
\newcommand{\eeq}{\end{equation}}
\newcommand{\beqs}{\begin{eqnarray}}
\newcommand{\eeqs}{\end{eqnarray}}

\catcode`@=11
\@addtoreset{equation}{section}
\@addtoreset{equation}{subsection}
\def\theequation{\ifnum\value{section}=0 \arabic{equation}\ignorespaces
\else \ifnum\value{section}=-1 A.\arabic{equation}\ignorespaces
\else \ifnum\value{subsection}=0 \thesection.\arabic{equation}\ignorespaces
\else \thesection.\arabic{subsection}.\arabic{equation}\ignorespaces
                           \fi
                      \fi
                 \fi}
\catcode`@=12

\begin{document}
\def\thefootnote{\fnsymbol{footnote}}

\def\square {{\vcenter {\hrule height .2mm
\hbox{\vrule width .2mm height 2mm \kern 2mm
\vrule width .2mm} \hrule height .2mm}}}
\outer
\def\endproof{{\unskip\nobreak\hfil\penalty50
\hskip2em\hbox{}\nobreak\hfil$\square$
\parfillskip0pt\finalhyphendemerits=0\par}
{\ifdim\lastskip<\bigskipamount\removelastskip
\penalty-10\bigskip\fi}}   

\baselineskip 6.0mm

\vspace{4mm}

\begin{center}

{\Large \bf Upper Bound for the Coefficients of Chromatic polynomials} 
\vspace{8mm}

\setcounter{footnote}{0}
Shu-Chiuan Chang 

\vspace{6mm}

C. N. Yang Institute for Theoretical Physics  \\
State University of New York       \\
Stony Brook, N. Y. 11794-3840  \\

\end{center}

\vspace{10mm}

proposed running head: Upper Bound for the Coefficients

email: shu-chiuan.chang@sunysb.edu

telephone number: 631-632-7982

fax number: 631-632-7954

\newpage

\begin{center}
{\bf Abstract}
\end{center}

This paper describes an improvement in the upper bound for the magnitude
of a coefficient of a term in the chromatic polynomial of a general graph.
If $a_r$ is the coefficient of the $q^r$ term in the chromatic polynomial
$P(G,q)$, where $q$ is the number of colors, then we find $a_r \le {e
\choose v-r} - {e-g+2 \choose v-r-g+2} + {e-k_g-g+2 \choose v-r-g+2} -
\sum _{n=1}^{k_g-\ell_g}\sum _{m=1}^{\ell_g-1} {e-g+1-n-m \choose v-r-g} -
\delta_{g,3}\sum _{n=1}^{k_g+\ell_{g+1}^*-\ell_g} {e-\ell_g-g+1-n \choose
v-r-g}$, where $k_g$ is the number of circuits of length $g$ and $\ell_g$
and $\ell_{g+1}^*$ are certain numbers defined in the text.

\vspace{10mm}

key words: chromatic polynomial

\vspace{16mm}

\pagestyle{empty}
\newpage

\pagestyle{plain}
\pagenumbering{arabic}
\renewcommand{\thefootnote}{\arabic{footnote}}
\setcounter{footnote}{0}

\section{Introduction} 

Let $G$ be a loopless graph with $v$ vertices and $e$ edges.  The
chromatic polynomial $P(G,q)$ counts the number of ways of coloring the
vertices of $G$ with $q$ colors subject to the condition that adjacent
vertices have different colors \cite{bbook, birk}.  (More generally, one
may consider multigraphs $G$ with multiple edges; however, it is
elementary that the chromatic polynomial for a graph with multiple edges
joining two vertices $v$ and $v^\prime$ is the same as for the graph with
just one edge joining these vertices.)  Besides its role in graph theory,
this polynomial is of interest in statistical physics as the
zero-temperature value of a certain model of cooperative phenomena and
phase transitions known as the Potts model (e.g., \cite{s4, wurev} and
references therein.).  The chromatic polynomial of a graph can be
calculated by means of the iterative use of the deletion-contraction
theorem or equivalently, the addition-contraction theorem, which
explicitly shows that it is a polynomial of maximal degree $v$ in $q$.  It
can be written as 
\beq 
P(G,q)=\sum _{r=1}^v (-1)^{v-r}a_rq^r \ .
\label{chrom}
\eeq
Since the chromatic polynomial of a set of graphs is the product of the
chromatic polynomials of each individual graph, we shall restrict our 
attention to a connected graph. 

\section{Basic Properties of $a_r$}

Since the calculation of $P(G,q)$ is, in general, a \#P-hard problem 
\cite{welsh}), it is useful to have bounds on the coefficients. The 
coefficients of a chromatic polynomial can be expressed as the sum of the
number of its subgraphs which do not contain any broken circuit 
\cite{boll, rrev}.

By an elementary application of the deletion-contraction theorem, it
follows that the coefficients $a_r$ in eq. (\ref{chrom}) are positive.
Furthermore, the leading terms are \cite{meredith}

\beq
a_r = \cases{{e \choose v-r} & if \quad $r>v-g+1$ \cr
             {e \choose v-r}-k_g & if \quad $r=v-g+1$\cr}
\label{leadingar}
\eeq
where $g$ is the girth of the graph, and $k_g$ is the number of circuits
of length $g$ in the graph. In particular, $a_v=1$ and $a_{v-1}=e$.

In practice, one finds that the $a_r$ increase monotonically as $r$
decreases from $v$, with at most two of these coefficients having the
maximal value, and then the $a_r$ decrease monotonically for lower values
of $r$.  We can easily show that the magnitude of $a_j$ can not just
increase monotonically without decreasing as follows,

\bigskip

{\bf Proposition 1} \quad The statement that $a_r \le a_{r-1} \quad 
\forall \ 1 \le r \le v$ is false except for the trivial case $e=1$.

{\sl Proof} \quad This statement is equivalent to 
\beq
a_v\le a_{v-1}\le a_{v-2} \le ... \le a_1 \ .
\label{monoar}
\eeq

Here we consider the graph with $e \ge 2$, and therefore, {\it a
fortiori}, $P(G,q=1)=0$. However, eq. (\ref{monoar}) implies
\beqs
P(G,q=1) & = & a_v - a_{v-1} + a_{v-2} - a_{v-3} + ... + (-1)^{v-1}a_1
\cr\cr & = & 1 - e + a_{v-2} - a_{v-3} + ... + (-1)^{v-1}a_1 \ne 0 \ .
\label{Pq1}
\eeqs
This contradiction disproves the statement except for the trivial case:
the tree graph with only two vertices and one edge has 
$P(T_2,q)=q(q-1)=q^2-q$.
\endproof

Recall the deletion-contraction theorem \cite{birk}-\cite{bbook}: Let $x$
and $y$ be adjacent vertices in $G$, and denote the edge joining them as
$xy$. Then

\beq
P(G,q)=P(G-xy,q)-P(G/xy,q)
\label{dc}
\eeq
where $G-xy$ is the graph obtained from $G$ by deleting the edge $xy$, and 
$G/xy$ is the graph obtained from $G$ by deleting the edge $xy$ and 
identifying $x$ and $y$.

   If we also write 
\beq
P(G-xy,q)=\sum _{r=1}^v (-1)^{v-r}a'_rq^r
\label{chromd}
\eeq
and
\beq
P(G/xy,q)=\sum _{r=1}^{v-1}(-1)^{v-1-r}a''_rq^r,
\label{chromc}
\eeq
then
\beq
a_r = a'_r + a''_r \quad {\rm for} \ \ 1 \le r \le v-1 \ .
\label{aradd}
\eeq

\section{Upper Bound on $a_r$}

It was known that the coefficients $a_r$'s are bounded above by the
corresponding coefficients of the complete graph with the same number of
vertices, $K_v$ \cite{rrev}. However, this upper bound is sharp only for
complete graphs. 

An upper bound on $a_r$ was given by Li and Tian \cite{litian} and is
\beq
a_r \le {e \choose v-r} - {e-g+2 \choose v-r-g+2} + {e-k_g-g+2 \choose
v-r-g+2} \ .
\label{ltub}
\eeq
We improve this bound. Let us use the convention
\beq
{a \choose b} = \cases{1 & if \quad $b=0$\cr
                       0 & if \quad $b>a$ or $b<0$\cr}
\label{choosecon}
\eeq
and, for some function $f(n)$
\beq
\sum _{n=0}^m f(n) = 0 \quad {\rm if} \quad m<0 \ .
\label{sumcon}
\eeq
Therefore, the bound reduced to the exact values in eq. (\ref{leadingar})
when $r \ge v-g+1$.

Let us derive a basic relation which will be used repeatedly later:

\medskip

{\bf Lemma 1} \quad If $a > b \ge c \ge 0$, then 
\beq
- {a \choose c} + {b \choose c} = -\sum _{n=1}^{a-b}{a-n \choose c-1} \ .
\label{lemma1}
\eeq

{\sl Proof} \quad We know
\beqs
{b+1 \choose c} - {b \choose c-1} & = & \frac{(b+1)!}{c!(b+1-c)!} -
       \frac{b!}{(c-1)!(b-c+1)!}\cr\cr
 & = & \frac{b!}{(c-1)!(b+1-c)!}\Bigl [\frac{b+1}{c} - 1\Bigr ]\cr\cr
 & = & \frac{b!}{(c-1)!(b+1-c)!}\frac{b+1-c}c\cr\cr
 & = & \frac{b!}{c!(b-c)!}\cr\cr
 & = & {b \choose c},
\label{lamma1proof}
\eeqs
and the result follows if we keep on applying this relation on the 
positive term generated from ${b\choose c}$.
\endproof

Next we prove a lemma that will be used for our bound.  To begin, we make
a choice of a certain edge $xy$ in $G$ where we shall apply the 
deletion-contraction theorem.  We then define $\ell_g$ as the number of 
circuits in $G$ of length $g$ that contain this edge $xy$. 
\medskip

{\bf Lemma 2} \quad If the number of circuits of length $n$ in a graph $G$
is $k_n$, where $k_n \ge 0$, $g \le n \le s$, $g \le s \le v$, and the
number of circuits of length $n$ containing the edge $xy$ is $\ell_n$,
then there are $v$ vertices and $e-1$ edges in graph $G-xy$, and the
number of circuits of length $n$ is $k'_n=k_n-\ell_n$. For the graph
$G/xy$, there are $v-1$ vertices, and the number of edges and the number
of circuits of length $n$ are (i) $e-1$ and $k''_n=k_n-\ell_n+\ell_{n+1}$
if $\ell_3=0$, or (ii) $e-\ell_3-1$ and $k''_n=k_n-\ell_n+\ell_{n+1}^*$ if
$\ell_3 \ne 0$, where $\ell_{n+1}^*$ is the number of circuits of length
$n+1$ which do not contain the edge $xz$ (or the edge $yz$) for any vertex
$z$.

{\sl Proof} \quad The number of vertices and edges is clear for $G-xy$,
and $G/xy$ when $\ell_3=0$. If $\ell_3\ne 0$, the contraction of the edge
$xy$ will result in $\ell_3$ double edges, and one of these edges can be
removed from each pair without affecting the chromatic polynomial. The
circuits of length $n+1$ in $G$ which become the circuits of length $n$ in
$G/xy$ and contain both edges $xy$ and $xz$ (or $yz$) are double-counted.
They are the same as the circuits of length $n$ in $G$ which contain the
edge $yz$ (or $xz$) but not the edge $xy$.  
\endproof

Now the upper bound of Li and Tian can be improved with extra negative
terms.

\bigskip

{\bf Theorem 1} \quad If the girth of a graph $G$ is $g$, and the number
of circuits of length $g$ in the graph is $k_g$, then
\beqs
& & a_r \le {e \choose v-r} - {e-g+2 \choose v-r-g+2} + 
            {e-k_g-g+2 \choose v-r-g+2}\cr\cr 
    & &-\sum _{n=1}^{k_g-\ell_g}\sum _{m=1}^{\ell_g-1}{e-g+1-n-m \choose v-r-g}
 -\delta_{g,3}\sum _{n=1}^{k_g+\ell_{g+1}^*-\ell_g} {e-\ell_g-g+1-n 
\choose v-r-g},
\label{arub}
\eeqs
where, as defined above, $\ell_g$ and $\ell_{g+1}^*$ are determined by the
initial choice of the edge $xy$ on which we apply the deletion-contraction 
theorem. 

{\sl Proof} \quad Consider the special case $g=3$ first. By eq. 
(\ref{ltub}) and Lemma 2
\beqs
a'_r & \le & {e-1 \choose v-r} - {(e-1)-3+2 \choose v-r-3+2} +
           {(e-1)-(k_3-\ell_3)-3+2 \choose v-r-3+2}\cr\cr
     & = & {e-1 \choose v-r} - {e-2 \choose v-r-1} +
           {e-2-k_3+\ell_3 \choose v-r-1}
\label{proofg3arp}
\eeqs
\beqs
a''_r & \le & {e-\ell_3-1\choose (v-1)-r}-{(e-\ell_3-1)-3+2 \choose 
(v-1)-r-3+2}
 + {(e-\ell_3-1)-(k_3-\ell_3+\ell_4^*)-3+2 \choose (v-1)-r-3+2}\cr\cr
      & = & {e-\ell_3-1 \choose v-r-1} - {e-\ell_3-2 \choose v-r-2} + 
            {e-k_3-\ell_4^*-2 \choose v-r-2}
\label{proofg3arpp}
\eeqs
then by eq. (\ref{aradd}) and Lemma 1,
\beqs
a_r & \le & {e-1 \choose v-r} - {e-2 \choose v-r-1 } +
            {e-2-k_3+\ell_3 \choose v-r-1}\cr\cr 
    & & + {e-\ell_3-1 \choose v-r-1} - {e-\ell_3-2 \choose v-r-2}
        + {e-k_3-\ell_4^*-2 \choose v-r-2}\cr\cr
    & = & {e \choose v-r} - {e-1 \choose v-r-1} - 
          \sum _{n=1}^{k_3-\ell_3} {e-2-n \choose v-r-2} +
          {e-\ell_3-1 \choose v-r-1}\cr\cr 
    & & - {e-k_3-1 \choose v-r-1} + {e-k_3-1 \choose v-r-1} -
      \sum _{n=1}^{k_3+\ell_4^*-\ell_3} {e-\ell_3-2-n \choose v-r-3}\cr\cr
    & = & {e \choose v-r} - {e-1 \choose v-r-1} + {e-k_3-1 \choose v-r-1}
          - \sum _{n=1}^{k_3-\ell_3} {e-2-n \choose v-r-2}\cr\cr
    & & + \sum _{n=1}^{k_3-\ell_3} {e-\ell_3-1-n \choose v-r-2} -
      \sum _{n=1}^{k_3+\ell_4^*-\ell_3} {e-\ell_3-2-n \choose v-r-3}\cr\cr
    & = &{e \choose v-r}-{e-1 \choose v-r-1}+{e-k_3-1 \choose v-r-1}\cr\cr
 & & - \sum _{n=1}^{k_3-\ell_3}\sum _{m=1}^{\ell_3-1}{e-2-n-m \choose v-r-3}
        - \sum _{n=1}^{k_3+\ell_4^*-\ell_3} {e-\ell_3-2-n \choose v-r-3} \ .
\label{proofg3ar}
\eeqs

   Consider $g>3$ and choose the edge $xy$ so that the girth of $G/xy$ is
$g-1$ and the number of circuits of length $g-1$ is $\ell_g$. By eq. 
(\ref{ltub}) and Lemma 2
\beqs
a'_r & \le & {e-1 \choose v-r} - {(e-1)-g+2 \choose v-r-g+2} +
             {(e-1)-(k_g-\ell_g)-g+2 \choose v-r-g+2}\cr\cr
     & = & {e-1 \choose v-r} - {e-g+1 \choose v-r-g+2} +
           {e-k_g+\ell_g-g+1 \choose v-r-g+2}
\label{proofarp}
\eeqs
\beqs
a''_r & \le & {e-1 \choose (v-1)-r}-{(e-1)-(g-1)+2\choose (v-1)-r-(g-1)+2}
              + {(e-1)-\ell_g-(g-1)+2 \choose (v-1)-r-(g-1)+2}\cr\cr
      & = & {e-1 \choose v-r-1} - {e-g+2 \choose v-r-g+2} +
            {e-\ell_g-g+2 \choose v-r-g+2}
\label{proofarpp}
\eeqs
therefore,
\beqs
a_r & \le & {e-1 \choose v-r} - {e-g+1 \choose v-r-g+2} +
            {e-k_g+\ell_g-g+1 \choose v-r-g+2}\cr\cr 
    & & + {e-1 \choose v-r-1} - {e-g+2 \choose v-r-g+2} + 
          {e-\ell_g-g+2 \choose v-r-g+2}\cr\cr
    & = & {e \choose v-r} - {e-g+2 \choose v-r-g+2} - 
          \sum _{n=1}^{k_g-\ell_g} {e-g+1-n \choose v-r-g+1}\cr\cr
    & & + {e-k_g-g+2 \choose v-r-g+2} - {e-k_g-g+2 \choose v-r-g+2} + 
          {e-\ell_g-g+2 \choose v-r-g+2}\cr\cr
    & = & {e \choose v-r} - {e-g+2 \choose v-r-g+2} + 
          {e-k_g-g+2 \choose v-r-g+2}\cr\cr
    & & - \sum _{n=1}^{k_g-\ell_g} {e-g+1-n \choose v-r-g+1}
        + \sum _{n=1}^{k_g-\ell_g} {e-\ell_g-g+2-n \choose v-r-g+1}\cr\cr
    & = & {e \choose v-r} - {e-g+2 \choose v-r-g+2} +
          {e-k_g-g+2 \choose v-r-g+2}\cr\cr 
  & &-\sum _{n=1}^{k_g-\ell_g}\sum _{m=1}^{\ell_g-1} {e-g+1-n-m \choose v-r-g} 
\ .
\label{proofar}
\eeqs
\endproof

It is obvious that this bound of $a_r$ is reduced to the Li-Tian bound in
eq. (\ref{ltub}) when $k_g=\ell_g=1$, and is optimized if we choose the
edge $xy$ so that the magnitude of the summation $S$, where 
\beq
S =  \sum _{n=1}^{k_g-\ell_g}\sum _{m=1}^{\ell_g-1} {e-g+1-n-m \choose v-r-g}
\label{sdef}
\eeq
is as large as possible. By Lemma 1, we can rewrite the bound as 
\beqs
a_r & \le & {e \choose v-r} - {e-g+2 \choose v-r-g+2} +
            {e-\ell_g-g+2 \choose v-r-g+2} - {e-g+1 \choose v-r-g+2}\cr\cr 
 & & + {e-k_g+\ell_g-g+1 \choose v-r-g+2} -\delta_{g,3}\Bigl [ {e-\ell_g-g+1
 \choose v-r-g+1} \cr\cr & & 
- {e-k_g-\ell_{g+1}^*-g+1 \choose v-r-g+1} \Bigr ] \ .
\label{arubn}
\eeqs

\vspace{10mm}

Acknowledgment: I would like to thank Prof. R. Shrock for helpful 
discussions.

\vfill
\eject
\end{document}